\documentclass[a4paper]{amsart}
\usepackage{latexsym,amssymb}
\usepackage[english]{babel}
\usepackage{bm}
\usepackage{amsfonts, amsthm}
\usepackage{subfig}
\usepackage{amsmath}
\usepackage{mathrsfs}
\usepackage{multirow}
\usepackage{bm}
\usepackage{pgf,tikz}
\usepackage{rotating}
\usepackage{calculator}

\usepackage{xcolor}

\newcommand{\comment}[1]{}

\usepackage{amsthm}
\theoremstyle{plain}

\newtheorem*{Theorem*}{Theorem 1}

\newtheorem*{Lemma*}{Lemma}
\newtheorem{Definition}{Definition}[section]
\newtheorem{Remark}[Definition]{Remark}
\newtheorem{Theorem}[Definition]{Theorem}

\newtheorem{Proposition}[Definition]{Proposition}

\newtheorem{Lemma}[Definition]{Lemma}
\newtheorem{Problem}[Definition]{Problem}
\newtheorem{Corollary}[Definition]{Corollary}

\title{
Variations on the Erd\H{o}s distinct-sums problem
}
\author{Simone Costa}
\address{DICATAM, Universit\`a degli Studi di Brescia, Via
Branze 43, 25123 Brescia, Italy}
\email{simone.costa@unibs.it}

\author{Marco Dalai}
\address{DII, Universit\`a degli Studi di Brescia, Via
Branze 38, 25123 Brescia, Italy}
\email{marco.dalai@unibs.it}

\author{Stefano Della Fiore}
\address{DII, Universit\`a degli Studi di Brescia, Via
Branze 38, 25123 Brescia, Italy}
\email{s.dellafiore001@unibs.it}
\begin{document}
\begin{abstract}
Let $\{a_1, . . . , a_n\}$ be a set of positive integers with $a_1 < \dots < a_n$ such that all $2^n$ subset sums are distinct. A famous conjecture by Erd\H{o}s states that $a_n>c\cdot 2^n$ for some constant $c$, while the best result known to date is of the form $a_n>c\cdot 2^n/\sqrt{n}$. In this paper, inspired by an information-theoretic interpretation, we extend the study to vector-valued elements $a_i\in \mathbb{Z}^k$ and we weaken the condition by requiring that only sums corresponding to subsets of size smaller than or equal to $\lambda n$ be distinct. For this case, we derive lower and upper bounds on the smallest possible value of $a_n$.
\end{abstract}

\maketitle

\section{Introduction}
For any $n\geq 1$, consider sets $\{a_1, . . . , a_n\}$ of positive integers with $a_1 < \dots < a_n$ whose subset sums are all distinct. A famous conjecture, due to Paul Erd\H{o}s, is that $a_n \geq c \cdot 2^n$ for some constant $c > 0$. Using the variance method, Erd\H{o}s and Moser \cite{Erdos} (see also \cite{Alon2}) were able to prove that
$$a_n\geq 1/4\cdot n^{-1/2}\cdot 2^n.$$
No advances have been made so far in removing the term $n^{-1/2}$ from this lower bound, but there have been several improvements on the constant factor, including the
work of Dubroff, Fox and Xu \cite{Fox}, Guy \cite{Guy}, Elkies \cite{Elkies}, Bae \cite{Bae}, and Aliev \cite{Aliev}. In particular, the best currently known lower bound states that
$$
a_n\geq (1+o(1))\sqrt{\frac{2}{\pi}} \frac{1}{\sqrt{n}}2^{n}\,.
$$
Two simple proofs of this result, first obtained unpublished by Elkies and Gleason, are presented in \cite{Fox}. In the other direction, the best-known
construction is due to Bohman \cite{Bohman}, who showed that there exist arbitrarily large such sets with $a_n \leq 0.22002 \cdot 2^n$.

In this paper, we propose a generalization of the problem in two directions. One is that the distinct-sums condition is weakened by only requiring that the sums of up to $\lambda n$ elements of the set be distinct, a direction with connections with the recent study independently proposed in \cite{MYR}. The second is that the integers $a_i$ be replaced by elements in $\mathbb{Z}^k$ for some $k\geq 1$. For these cases we derive both upper and lower bounds on the smallest possible value of the largest component among all of the $a_i$'s, that is on the smallest cube which contains all the $a_i$ elements.

This variation on the problem is inspired by an information-theoretic interpretation, namely in the setting of signaling over a multiple access channel. Looking at the original problem, we can interpret the $a_i$ integers as pulse amplitudes that $n$ transmitters can transmit over an additive channel to send one bit of information each, for example, to signal to the base station that they want to start a communication session. The requirement that all subset sums be distinct expresses the desire that the base station be able to infer any possible subset of active users. In this setting, a natural assumption to consider is that only a maximum fraction of the users might actually be active at the same time, and that signals be vector-valued rather than scalars since the channel would be used over an interval of time sending a sequence of pulses (codewords) rather than a single pulse.


More formally, we consider the following problem.
\begin{Problem}\label{SmallSize}
Let $\mathcal{F}_{\lambda,n}$ be the family of all subsets of $\{1,\dots,n\}$ whose size is smaller than or equal to $\lambda n$. We are interested in the minimum $M$ such that there exists a sequence $\Sigma=(a_1,\dots,a_n)$ in $\mathbb{Z}^k$, $a_i\in [0,M]^k\ \forall i$, (i.e. $\Sigma$ is $M$-bounded) such that for all distinct $A_1,A_2\in\mathcal{F}_{\lambda,n}$, $S(A_1)\neq S(A_2)$, where
$$S(A)=\sum_{i\in A} a_i\,.$$
In the following, we will call such sequences $\mathcal{F}_{\lambda,n}$-sum distinct.
\end{Problem}

Throughout the paper,  the logarithms are in base two and we denote the open interval with endpoints $x$ and $y$ by $(x,y)$ and the closed interval by $[x,y]$.

The paper is organized as follows.
Section $2$ is devoted to lower bounds on the values of $M$ in Problem \ref{SmallSize}. We show that for $\lambda\geq 1/2$, both the isoperimetric approach (see \cite{Fox}) and the variance method can be applied to obtain non-trivial lower bounds. Then, in Section $3$, we derive three upper bounds using, respectively, the combinatorial nullstellensatz, the probabilistic method, and a direct construction. \section{Lower Bounds}
In this section we will derive three different lower bounds on $M$. Firstly, we provide a very elementary (but still interesting since, for $\lambda<1/2$ we have no better results) lower bound.

\begin{Proposition}\label{prop:trivialBound}
Let $\Sigma=(a_1,\dots,a_n)$ be an
$\mathcal{F}_{\lambda,n}$-sum distinct sequence in $\mathbb{Z}^k$ that is $M$-bounded. Then
$$
M \geq (1 + o(1)) \cdot \begin{cases} \frac{1}{\lceil \lambda n \rceil \sqrt[k]{2 \pi n \lambda (1-\lambda)}} 2^{n h(\lambda)/k} & \mbox{ if } \lambda < 1/2;\\
\frac{1}{\lceil \lambda n \rceil} \cdot 2^{(n-1)/k} &\mbox{ if } 1/2 \leq \lambda < 1; \\
\frac{1}{n} \cdot 2^{n/k} &\mbox{ if } \lambda = 1;
\end{cases}\,
$$
where $h(\lambda) = -\lambda \log \lambda - (1-\lambda) \log (1-\lambda)$ is the binary entropy function. 
\end{Proposition}
\proof
The maximum possible sum we can get on some coordinates is at most $\lceil \lambda n \rceil M$. Then by the pigeonhole principle, for values of $\lambda \in (0,1/2)$, we have that
\begin{equation*}\label{eq:sumbinomial}
M^k \geq \frac{1}{\lceil \lambda n \rceil^k} \sum_{i=0}^{\lceil \lambda n \rceil} \binom{n}{i} \geq \frac{1}{\lceil \lambda n \rceil^k \sqrt{2 \pi n \lambda (1-\lambda)}} 2^{n h(\lambda)/k}\,.
\end{equation*}
This leads to the asymptotic bound as $n\to\infty$
\begin{equation*}\label{eq:trivallowerbound}
M \geq (1 + o(1)) \frac{1}{\lceil \lambda n \rceil \sqrt[k]{2 \pi n \lambda (1-\lambda)}} 2^{n h(\lambda)/k}.
\end{equation*}

For values of $\lambda \in [1/2, 1]$ the lower bound on $M$ can be easily derived noticing that the sum $\sum_{i=0}^{\lceil \lambda n \rceil} \binom{n}{i}$ is greater than or equal to $2^{n-1}$ for $\lambda \in [1/2, 1)$ and it is equal to $2^n$ for $\lambda = 1$. Therefore, we have that
\begin{equation}\label{eq:trivialn2}
M \geq (1 + o(1)) \cdot \begin{cases}
\frac{1}{\lceil \lambda n \rceil} \cdot 2^{(n-1)/k} &\mbox{ if } 1/2 \leq \lambda < 1; \\
\frac{1}{n} \cdot 2^{n/k} &\mbox{ if } \lambda = 1.
\end{cases}\,
\end{equation}
\endproof

Now, if $\lambda \geq 1/2$, we see that it is possible to improve on the term $C_n = 1 / \lceil \lambda n \rceil$ in \eqref{eq:trivialn2} using the Harper isoperimetric inequality (see \cite{Harper}) as done in \cite{Fox} for $\lambda=1$. In particular, we see that the same bound obtained for $\lambda=1$ also holds for all $\lambda>1/2$. For $\lambda=1/2$, instead, a weakening by a factor of $2$ appears, which can be explained in terms of the concentration of measure around the average value of the sums.

\begin{Theorem}\label{harper}[Harper vertex-isoperimetric inequality]
Let $\mathcal{G}$ be a family of subsets of $[1,n]$ with cardinality $\sum_{i=0}^k \binom{n}{i} \leq |\mathcal{G}| \leq 2^{n-1}$ then $|\partial \mathcal{G}| \geq \binom{n}{k+1}$ where ${\partial \mathcal{G} = \{ F \: | \: F \in \mathcal{P}([1,n]), \min_{Y \in \mathcal{G}} |F \Delta Y| = 1 \}}$ is called the border of $\mathcal{G}$.
\end{Theorem}

Inspired by \cite{Fox}, we obtain the following theorem.
\begin{Theorem}\label{harperbound}
Let $\Sigma=(a_1,\dots,a_n)$ be an
$\mathcal{F}_{\lambda,n}$-sum distinct sequence in $\mathbb{Z}$ that is $M$-bounded. Then
\begin{equation*}\label{eq:harpern2}
M \geq (1 + o(1)) \cdot \begin{cases} \frac{1}{\sqrt{2\pi n}} \cdot 2^n & \text{if } \lambda=1/2; \\ \sqrt{\frac{2}{\pi n}} \cdot 2^n &\text{if } \lambda \in (1/2, 1]. \end{cases}
\end{equation*}
\end{Theorem}
\proof
Assume that there exists an
$\mathcal{F}_{\lambda,n}$-sum distinct sequence $\Sigma=(a_1, a_2, \ldots, a_n)$ and, without loss of generality, that $a_1<a_2<\dots<a_n$.
Let $\mathcal{G}$ be a set of vectors $\epsilon = (\epsilon_1, \ldots, \epsilon_n)$ such that $\epsilon_i \in \{-1/2, 1/2\}$ and the dot product $\epsilon \cdot \Sigma <0\ \forall \epsilon\in \mathcal{G}$. Clearly $|\mathcal{G}| = 2^{n-1}$ by symmetry. Then by Theorem \ref{harper} we know that $|\partial \mathcal{G}| \geq \binom{n}{\lceil n/2 \rceil}$. If we take $\eta \in \partial \mathcal{G}$ then $0< \eta \cdot \Sigma < a_n$. We can express $\partial \mathcal{G} = \partial \mathcal{G}_1 \cup \partial \mathcal{G}_2$ where
$$\partial \mathcal{G}_1= \{ \eta \in\partial \mathcal{G} : \text{supp}(\eta +1/2) \leq \lfloor \lambda n \rfloor \}$$
and
$$\partial \mathcal{G}_2=\{ \eta \in \partial \mathcal{G} : \text{supp}(\eta +1/2) \geq \lfloor \lambda n \rfloor + 1 \}.$$

If $\lambda \in (1/2, 1]$, then we have that
\begin{equation}\label{eq:deltag1}
|\partial \mathcal{G}_1| \geq \binom{n}{\lceil n/2 \rceil} - |\partial \mathcal{G}_2|.
\end{equation}
Because of the definition of $\partial \mathcal{G}_2$
$$
|\partial \mathcal{G}_2| \leq \sum_{i=\lfloor \lambda n \rfloor+1}^n \binom{n}{i} \leq 2^{h(\lambda) n}.
$$
Since in this case $h(\lambda) < 1$, from \eqref{eq:deltag1} we obtain
$$
|\partial \mathcal{G}_1| \geq (1 + o(1))\binom{n}{\lceil n/2 \rceil}.
$$
Again, by the pigeonhole principle there exists $\eta_1, \eta_2 \in \partial \mathcal{G}_1$ such that
$$|(\eta_1 - \eta_2) \cdot \Sigma| <  a_n / |\partial \mathcal{G}_1| \leq (1 + o(1)) a_n / \binom{n}{\lceil n/2 \rceil}.$$
Finally, by the hypothesis of sum-distinctness we have that $|(\eta_1 - \eta_2) \cdot \Sigma| \geq 1$, and hence
$$ a_n>(1 + o(1))\binom{n}{\lceil n/2 \rceil}=(1 + o(1)) \sqrt{\frac{2}{\pi n}} \cdot 2^n.$$

For $\lambda = 1/2$ we need a tweak. In this case we see that either $\partial \mathcal{G}_1$ or $\partial \mathcal{G}_2$ is greater than or equal to $(1/2) {\binom{n}{\lceil n/2 \rceil}}$. Here we note that, since $\Sigma$ is $\mathcal{F}_{1/2,n}$-sum distinct, it is also $\overline{\mathcal{F}_{1/2,n}}$-sum distinct, where $\overline{\mathcal{F}_{1/2,n}}$ is the complement of $\mathcal{F}_{1/2,n}$ in the power set $\mathcal{P}([1,n])$. Therefore we can assume, without loss of generality, that $\partial \mathcal{G}_1$ is greater than or equal to $(1/2) {\binom{n}{\lceil n/2 \rceil}}$. Proceeding as in the previous case, here we obtain that
$$ a_n>(1/2)\binom{n}{\lceil n/2 \rceil}=(1 + o(1))\frac{1}{\sqrt{2\pi n}} \cdot 2^n.$$
\endproof
\begin{Remark}\label{rem:harper}
A simple extension of Theorem \ref{harperbound} to the case $k>1$ leads, for $\lambda > 1/2$, to the bound
$$
M \geq (1+o(1)) \sqrt[k]{\frac{2}{\pi}} n^{\frac{1}{2k}-1} 2^{n/k}.
$$
Although a more refined reasoning might lead to better results, we did not manage to obtain something which could compete with Theorem \ref{teo:multidimguy} below.
\end{Remark}
Now we see that, using the variance method (see \cite{Alon2}, \cite{Erdos} or \cite{Guy}), it is possible to improve the bound of Remark \ref{rem:harper} whenever $k>1$ and $\lambda \in [1/2, 1]$.
\begin{Theorem}\label{teo:multidimguy}
Let $\lambda \geq 1/2$ and let $\Sigma=(a_1,\dots,a_n)$ be an
$\mathcal{F}_{\lambda,n}$-sum distinct sequence in $\mathbb{Z}^k$ that is $M$-bounded. Then
$$
M \geq (1+o(1)) \cdot \begin{cases} \sqrt{\frac{4}{\pi n (k+2)}} \cdot \Gamma(k/2+1)^{1/k} \cdot 2^{n/k} & \text{if } \lambda=1; \\
\sqrt{\frac{4}{\pi n (k+2)}} \cdot \Gamma(k/2+1)^{1/k} \cdot 2^{(n-1)/k} & \text{if } 1/2 \leq \lambda < 1; \end{cases}\,
$$
where $\Gamma$ is the gamma function.
\end{Theorem}
\proof
Let $\Sigma=(a_1,\dots,a_n)$ be an $M$-bounded and $\mathcal{F}_{\lambda,n}$-sum distinct sequence in $\mathbb{Z}^k$ where $\lambda\geq 1/2$.
Consider a random variable $X = \sum_{i=1}^n \epsilon_i a_i$ where the random vectors $(\epsilon_1, \epsilon_2, \ldots, \epsilon_n)$ are uniformly distributed over the set $\{ \epsilon \in \{-1/2, 1/2\}^n : \text{supp}(\epsilon+1/2) \leq \lambda n \}$. We denote with $\mu$ and $\sigma^2$ respectively the expected value and the variance of the random variable $X$.

We know that $\sigma^2 = \mathbb{E}[|X|^2] - |\mathbb{E}[X]|^2 \leq \mathbb{E}[|X|^2]$. Expanding $\mathbb{E}[|X|^2]$ we get
\begin{equation}\label{eq:varianceX}
\mathbb{E}[|X|^2] = 1/4 \sum_{i=1}^n |a_i|^2 + 2 \sum_{i < j} \mathbb{E}[\epsilon_i \epsilon_j] (a_i \cdot a_j),
\end{equation}
where $\mathbb{E}[\epsilon_i \epsilon_j]$ does not depend on the specific values chosen for $i$ and $j$. Then for each $i \neq j$ the following inequality holds
\begin{align}\label{eq:crossProduct}
\mathbb{E}[\epsilon_i \epsilon_j] &= 1/4 \cdot \frac{\sum_{i=0}^{\lfloor \lambda n \rfloor} \binom{n-2}{i} + \sum_{i=0}^{\lfloor \lambda n \rfloor - 2} \binom{n-2}{i} - 2 \sum_{i=0}^{\lfloor \lambda n \rfloor -1} \binom{n-2}{i}}{\sum_{i=0}^{\lfloor \lambda n \rfloor} \binom{n}{i}} \nonumber \\
& = 1/4 \cdot \frac{\binom{n-2}{ \lfloor \lambda n \rfloor} - \binom{n-2}{\lfloor \lambda n \rfloor - 1}}{\sum_{i=0}^{\lfloor \lambda n \rfloor} \binom{n}{i}} \nonumber \\
& \leq 0\,,
\end{align}
where the inequality holds since $\lambda \geq 1/2$.
By \eqref{eq:varianceX}, \eqref{eq:crossProduct}, since $|a_i|^2 \leq k M^2$, we have that
\begin{equation}\label{eq:upperSigma}
\sigma^2 \leq \mathbb{E}[|X|^2] \leq \frac{k n M^2}{4}\,.
\end{equation}

Now we want to provide a lower bound on $\sigma^2$. We know, by the sum-distinctness property of $\Sigma$, that each possible value of $X$, has a probability of happening equal to $1/|\mathcal{F}_{\lambda,n}|$. Therefore, considered the possible outcomes $s_1,s_2,\dots, s_{|\mathcal{F}_{\lambda,n}|}$ of the random variable $X$ and its mean $\mu$, the variance can be expressed as follows
\begin{align*}
\sigma^2 = \frac{1}{|\mathcal{F}_{\lambda,n}|} \sum_{i=1}^{|\mathcal{F}_{\lambda,n}|} |s_i - \mu|^2.
\end{align*}
Thus, we can lower bound the variance by the minimum value that the above expression can take for distinct values of the $s_i$'s on a discrete grid when we relax the constraint that $\mu$ be their average. For any fixed $\mu$, the sum is minimized when the $s_i$'s are packed as close as possible around $\mu$, that is, if $d=\max_{i}|s_i-\mu|$, then no point in the grid at a distance $d'<d$ from $\mu$ is left unused (otherwise we can move one of the $s_i$'s closer to $\mu$ and make the sum smaller). Let $R$ be the radius of a ball of volume $|\mathcal{F}_{\lambda,n}|$, that is,
\begin{equation}\label{eq:sphereRadius}
R = \frac{\Gamma(k/2 + 1)^{1/k}}{\sqrt{\pi}} |\mathcal{F}_{\lambda,n}|^{1/k}.
\end{equation}
By considering unit-volume non-overlapping cubes around each point in the grid, we deduce that $d\geq R'=R-\sqrt{k}$, so that we have an $s_i$ in any discrete point at distance $d'<R'$ from $\mu$. So, we have
\begin{align*}
\sigma^2 \geq \frac{1}{|\mathcal{F}_{\lambda,n}|} \sum_{|s-\mu|< R'}|s - \mu|^2
\end{align*}
where $s$ runs over all points in the ball on a discrete grid with spacing $1$. If we thus scale everything down by $R'$, renaming $\tilde{s}$ and $\tilde{\mu}$ the scaled quantities, we find
\begin{align*}
\sigma^2  & \geq  \frac{R'^2}{|\mathcal{F}_{\lambda,n}|} \sum_{|\tilde{s}-\tilde{\mu}|< 1}|\tilde{s} - \tilde{\mu}|^2\\
& =   \frac{R'^{2+k}}{|\mathcal{F}_{\lambda,n}|} \sum_{|\tilde{s}-\tilde{\mu}|< 1}|\tilde{s} - \tilde{\mu}|^2\frac{1}{R'^k}
\end{align*}
 where now $\tilde{s}$ runs over all points in the ball on a discrete grid with spacing $1/R'$. For fixed $k$, as $n\to\infty$  $R'$ grows to infinity with $R'=(1+o(1))R$, and the sum in the last expression behaves as a Riemann approximation for an integral over a unit ball. So, asymptotically as $n\to\infty$ we have 
$$
\sigma^2 \geq (1+o(1))\frac{R^{2+k}}{|\mathcal{F}_{\lambda,n}|} \int_{|\tilde{x}-\tilde{\mu}| \leq 1} |\tilde{x}-\tilde{\mu}|^2 \, d\tilde{x}\,.
$$
Integrating in polar coordinates, using the $(k-1)$-dimensional volume of the $(k-1)$-dimensional sphere of radius $\rho$, $S_{k-1}(\rho) = \frac{k\pi^{k/2}}{\Gamma(k/2 + 1)} \rho^{k-1}$, we obtain
\begin{align*}
\sigma^2 & \geq (1+o(1)) \frac{R^{2+k}}{|\mathcal{F}_{\lambda,n}|} \int_{0}^{1} S_{k-1}(\rho) \rho^2 d\rho\\
& \geq (1+o(1))\frac{R^{2+k}}{|\mathcal{F}_{\lambda,n}|} \frac{k\pi^{k/2}}{\Gamma(k/2 + 1)(k+2)}.
\end{align*}
Using \eqref{eq:sphereRadius} and \eqref{eq:upperSigma} we obtain the thesis.
\endproof
\begin{figure}%
\centering
\subfloat[\centering $k=1$]{{\includegraphics[width=5.8cm]{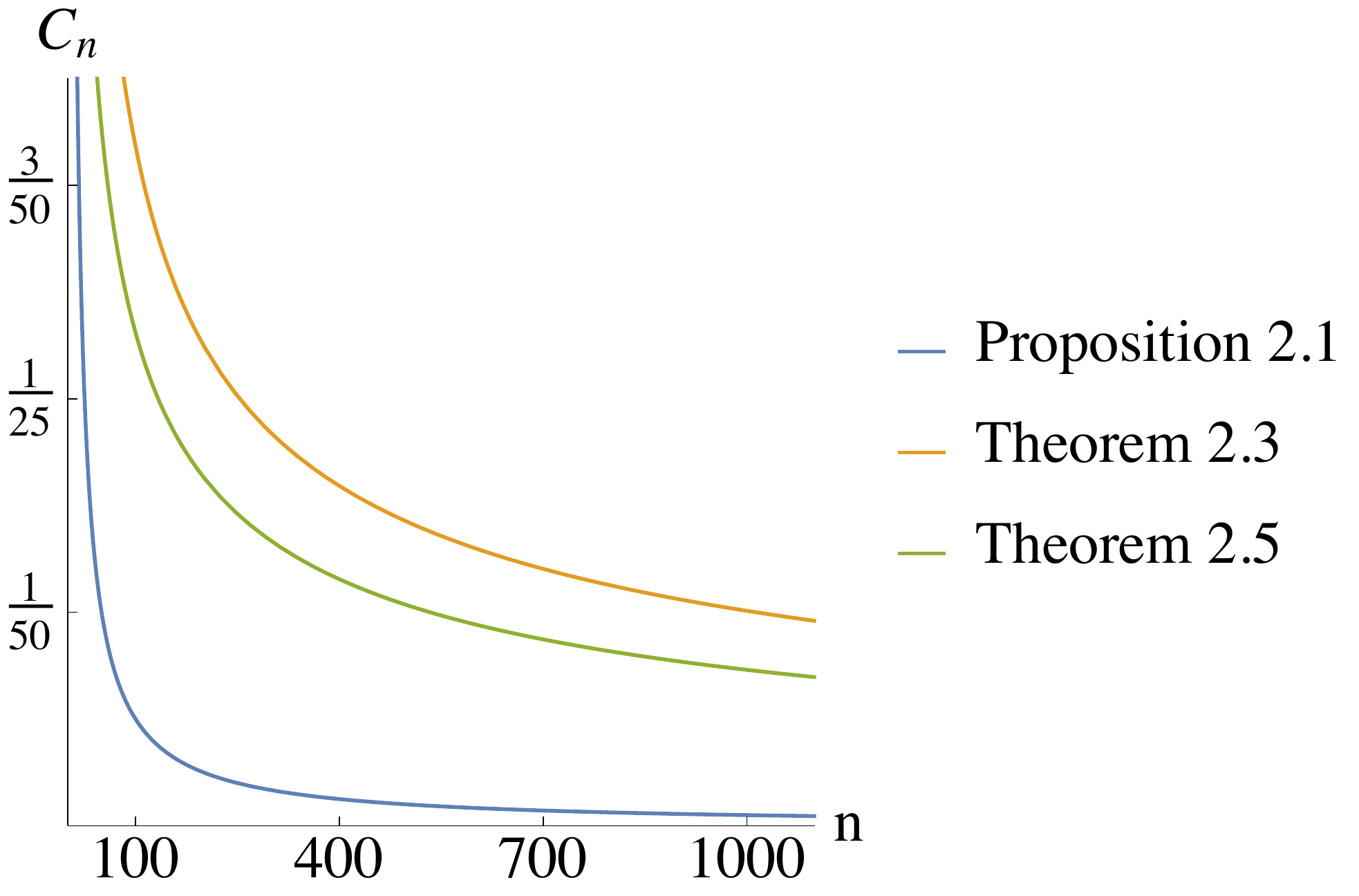} }}%
\qquad
\subfloat[\centering $k>1$]{{\includegraphics[width=5.8cm]{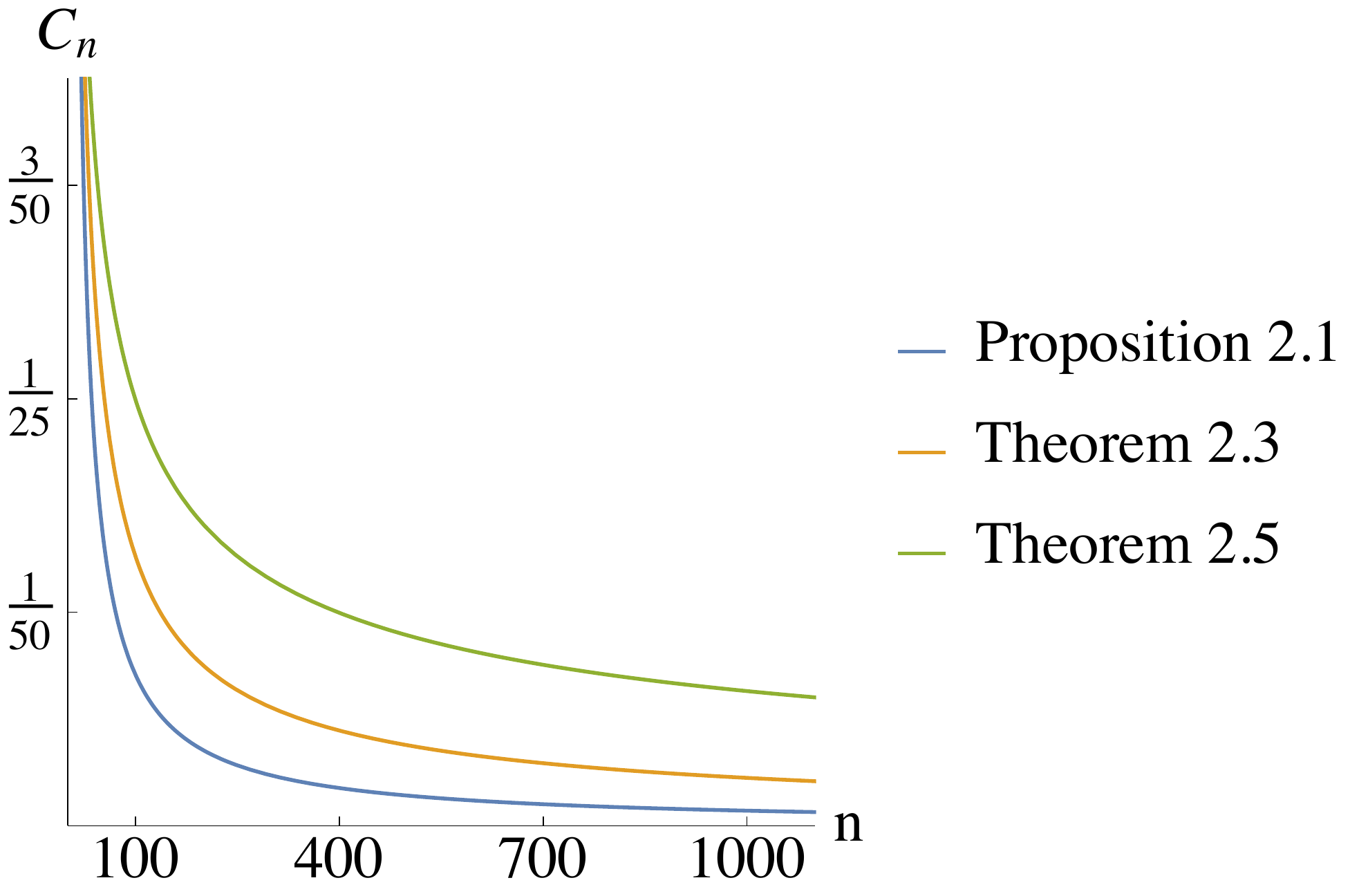} }}%
\caption{Representation of the sub-exponential factor $C_n$ of the lower bounds for $1/2 \leq \lambda \leq 1$.}%
\label{fig:example}%
\end{figure}

\section{Upper Bounds}
The goal of this section is to provide upper bounds on $M$. We remark that the best known upper bound for the classical Erd\H{o}s distinct-sums problem (see Bohman \cite{Bohman}) is always (i.e. for any $\lambda$) an upper bound on $M$ for the $\mathcal{F}_{\lambda,n}$ distinct-sums problem. Now we will see that this bound can be improved in several situations.
\subsection{One dimensional Upper Bounds}
In this paragraph, we consider the one-dimensional case that is $k=1$. In this case, we can provide an upper bound by using Alon's combinatorial nullstellensatz. This Theorem has been applied in several Combinatorial Number Theory problems; we refer to \cite{HOS} (see also \cite{CP}) for applications in the similar context of Alspach's partial sums conjecture and to \cite{CMPP} for background on that problem. We report here the theorem for the reader's convenience.

\begin{Theorem}
\cite[Theorem 1.2]{Alon}\label{thm:alon}
Let $\mathbb{F}$ be a field and let $f=f(x_1,\ldots,x_k)$ be a polynomial in $\mathbb{F}[x_1,\ldots,x_k]$. Suppose the degree of $f$ is
$\sum\limits_{i=1}^k t_i$,
where each $t_i$ is a nonnegative integer, and suppose the coefficient of
$\prod\limits_{i=1}^k x_i^{t_i}$ in $f$ is nonzero.
Then, if $A_1,\ldots,A_k$ are subsets of $\mathbb{F}$ with $|A_i|> t_i$, there are $a_1 \in A_1,\ldots,a_k \in A_k$ so that
$f(a_1,\ldots,a_k)\neq 0$.
\end{Theorem}

Before providing our upper bound, we need an enumerative lemma.
The bound that we will derive is non-trivial,  i.e.,  it is better than the one derived from the powers of two sequence,  for $\lambda < \bar{\lambda} \approx 0.113546$,  so we assume for simplicity that $\lambda<1/3$.

We define for convenience
$$
f(\lambda) = H(\lambda, \lambda, 1-2\lambda),
$$
where $H(p_1, \ldots, p_h) = \sum_{i=1}^h - p_i \log p_i$ is the Shannon entropy of a probability vector $(p_1, \ldots, p_h)$.

\begin{Lemma}\label{Stima2}
Let $\mathcal{C}_{\bar{i}}$ be the family of the unodered pairs $\{A_1,A_2\}$ of subsets of $[1,n]$ such that, given an element $\bar{i}\in [1,n]$:
\begin{itemize}
\item $A_1\cap A_2=\emptyset$;
\item The element $\bar{i}$ belongs to $A_1\cup A_2$;
\item The cardinalities of $A_1$ and $A_2$ are smaller than or equal to $\lambda n$.
\end{itemize}
Then, for $\lambda<1/3$, we have the following upper bound on the cardinality of $\mathcal{C}_{\bar{i}}$

\begin{equation*}\label{eq:multibound}
|\mathcal{C}_{\bar{i}}| < \lambda^3 n^2 \cdot 2^{ f(\lambda) n}\,.
\end{equation*}
\end{Lemma}
\proof
Suppose, without loss of generality, that $\bar{i} \in A_1$. Then we can upper bound the size of $\mathcal{C}_{\bar{i}}$ as follows
$$
|\mathcal{C}_{\bar{i}}| \leq \sum_{i=1}^{\lfloor\lambda n\rfloor} \sum_{j=0}^{\lfloor\lambda n\rfloor} \binom{n-1}{i-1, j, n-i-j}
$$
where $i$ represents the cardinality of $A_1$ while $j$ that of $A_2$. Using the fact that $\binom{n-1}{i-1, j, n-i-j} \leq \lambda \cdot \binom{n}{i, j, n-i-j}$ for each $i \in [1, \lambda n]$ and $j \in [0, \lambda n]$, we get
\begin{equation*}
|\mathcal{C}_{\bar{i}}| < \lambda^3 n^2 \binom{n}{\lfloor\lambda n\rfloor, \lfloor\lambda n\rfloor, n-2\lfloor\lambda n\rfloor}\,,
\end{equation*}
since the multinomial coefficient is maximized when all numbers are as equal as possibile.
Then, by a well-known entropy bound on the multinomial coefficient (see \cite[Lemma 2.2]{Csiszar}) we have that
\begin{equation*}
|\mathcal{C}_{\bar{i}}| < \lambda^3 n^2 \cdot 2^{n H(\frac{\lfloor\lambda n\rfloor}{n}, \frac{\lfloor\lambda n\rfloor}{n}, 1-2 \frac{\lfloor\lambda n\rfloor}{n})}\leq \lambda^3 n^2 \cdot 2^{n H(\lambda, \lambda, 1-2 \lambda)}\,
\end{equation*}
where the last inequality holds because, for $\lambda < 1/3$, $f(\lambda)$ is an increasing function.
\endproof


We are now ready to state our bound.

\begin{Theorem}\label{NUB}
For any $\lambda<1/3$, there exists a sequence $\Sigma=(a_1,\dots,a_n)$ of $\left(\lambda^3 n^2 2^{f(\lambda) n}\right )$-bounded positive integers that is $\mathcal{F}_{\lambda,n}$-sum distinct.
\end{Theorem}
\proof
For any pair $(A_1,A_2)\in \mathcal{F}_{\lambda,n}^2$, we define the linear polynomial $$l_{A_1,A_2}(x_1,\dots,x_n):=\sum_{i\in A_1} x_i-\sum_{j\in A_2} x_j.$$
Now, let us denote by $\mathcal{P}_{\lambda,n}$ the family of the pairs $(A_1,A_2)$ of elements of $\mathcal{F}_{\lambda,n}$ such that $A_1\cap A_2=\emptyset$ and $\min(A_1)<\min(A_2)$.
Then we set
$$q_{\mathcal{F}_{\lambda,n}}(x_1,\dots,x_n):=\prod_{(A_1,A_2)\in \mathcal{P}_{\lambda,n}}l_{A_1,A_2}(x_1,\dots,x_n).$$
We note that, for any pair $(A'_1,A'_2)\in \mathcal{F}_{\lambda,n}^2$ such that $A'_1\not=A'_2$, the linear polynomial $l_{A'_1,A'_2}(x_1,\dots,x_n)$ is equal to $\pm l_{A_1,A_2}(x_1,\dots,x_n)$ for some $(A_1,A_2)\in \mathcal{P}_{\lambda,n}$.
Therefore $\Sigma=(a_1,\dots,a_n)$ is $\mathcal{F}_{\lambda,n}$-sum distinct if and only if $q_{\mathcal{F}_{\lambda,n}}(a_1,\dots,a_n)\not=0$.

Since $\mathbb{Z}[x_1,\dots,x_n]$ is an integral domain, $q_{\mathcal{F}_{\lambda,n}}$ is not constantly zero. Therefore there exist $t_1,\dots,t_n$,
where each $t_i$ is a nonnegative integer, such that the coefficient of
$\prod\limits_{i=1}^n x_i^{t_i}$ in $q_{\mathcal{F}_{\lambda,n}}$ is nonzero.
Since $q_{\mathcal{F}_{\lambda,n}}$ is homogeneous, we also have that its degree is $\sum\limits_{i=1}^n t_i$. Let us consider $\bar{i}$ such that $t_{\bar{i}}=\max_i t_i$. The term $x_{\bar{i}}^{t_{\bar{i}}}$ originates from the factor $r_{\bar{i}}$ of $q_{\mathcal{F}_{\lambda,n}}$ defined by the product
$$r_{\bar{i}}(x_1,\dots,x_n):=\prod_{(A_1,A_2)\in \mathcal{P}_{\lambda,n}:\ \bar{i}\in A_1\cup A_2}l_{A_1,A_2}(x_1,\dots,x_n).$$
Hence, because of Lemma \ref{Stima2}, we have that $t_{\bar{i}}< \lambda^3 n^2 \cdot 2^{f(\lambda)n}$. This means that the hypotheses of Theorem \ref{thm:alon} are satisfied whenever $M\geq \lambda^3 n^2 \cdot 2^{f(\lambda) n}>\max_i t_i$ and hence, under this constraint, there exist $a_1\in [1,M],\ldots,a_n \in [1,M]$ such that
$q_{\mathcal{F}_{\lambda,n}}(a_1,\ldots,a_n)\neq 0$.
\endproof

We recall that the result of Theorem \ref{NUB} is non-trivial only when $\lambda < \bar{\lambda} \approx 0.113546$. Now, we investigate the range $\lambda\in [\bar{\lambda},1/4)$; here we provide a direct construction that improves the constant of Bohman \cite{Bohman} bound.

\begin{Lemma}\label{Insertion1}
Let us consider the sequence $\widetilde{\Sigma}_{n}=(b_1,\dots,b_n)$ where
$$b_i:=\begin{cases}
2^{i-1} & i=1,2\dots n-1\\
\sum_{j=0}^{j<n/2-1} 2^{2j} & i=n;
\end{cases}$$
Then, given two subsets $A_1,A_2$ of $[1,n]$ such that $|A_1|+|A_2|< n/2$,
$$S(A_1)=\sum_{i\in A_1} b_i\not=\sum_{j\in A_2} b_j=S(A_2).$$
\end{Lemma}

\begin{figure}
\centering
\includegraphics[scale=1]{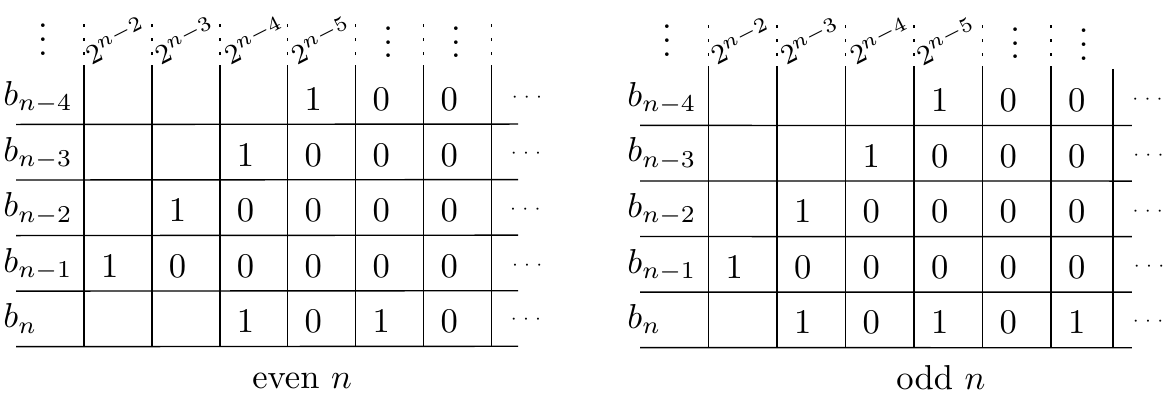}
\caption{Binary representation of the $b_n$ integers used in Lemma \ref{Insertion1}.}
\label{fig:table_b_n}
\end{figure}

The structure of the set $\widetilde{\Sigma}_{n}$ is better understood by writing a table of the binary representations of the integers $b_n$, as shown in Figure \ref{fig:table_b_n}.

\proof
Let us suppose, by contradiction, that there exist $n$, $A_1$ and $A_2$ with $|A_1|+|A_2|<n/2$ such that $S(A_1)=S(A_2)$, and let us consider the smallest $n$ for which this holds.

We note that if two sets $A_1$ and $A_2$ have the same sum, then also $A_1\setminus (A_1\cap A_2)$ and $A_2\setminus (A_1\cap A_2)$ have the same sum.
Therefore we may also assume that $A_1$ and $A_2$ are disjoint. Since a simple check shows that the thesis is true for $n \leq 5$, $n$ must be bigger than $5$.
Moreover, since $\widetilde{\Sigma}_{n}\setminus \{b_n\}$ is clearly sum-distinct, we can assume without loss of generality that $n\in A_1$. 
Therefore we have
$$b_n+\sum_{i\in A_1\setminus \{n\}} b_i=\sum_{j\in A_2}b_j.$$
which can be rewritten as
\begin{equation}
\sum_{i=0}^{ i<n/2-1} 2^{2i}+\sum_{i\in A_1\setminus \{n\}} 2^{i-1}=\sum_{j\in A_2}2^{j-1}.
\label{eq:general_n}
\end{equation}
Now we divide the proof in two cases, according to whether $n$ is even or $n$ is odd. The binary representations shown in Figure \ref{fig:table_b_n} might be useful as a complement in some steps of the discussion.

Consider the case of even $n$. First observe that in this case equation \eqref{eq:general_n} can rewritten by replacing $n$ with $n-1$ in the upper extreme of the first summation, that is,
\begin{equation}
\sum_{i=0}^{ i<(n-1)/2-1} 2^{2i}+\sum_{i\in A_1\setminus \{n\}} 2^{i-1}=\sum_{j\in A_2}2^{j-1}.
\label{eq:even_n}
\end{equation}
We now claim that $n-1\in A_2$. Indeed, if $n-1$ is neither in $A_1$ nor in $A_2$,
we see that equation \eqref{eq:even_n} provides a counterexample which is already contained in $\widetilde{\Sigma}_{n-1}$. Formally, the sets $A'_1=A_1\setminus \{n\}\cup \{n-1\}$ and $A'_2=A_2$ give a counterexample for $\widetilde{\Sigma}_{n-1}$ satisfying $|A'_1|+|A'_2|=|A_1|+|A_2|<(n-1)/2$, because $|A_1|+|A_2|<n/2$ with even $n$. This contradicts the minimality of $n$. It is easy to see that $n-1\in A_1$ is impossible, since we would have $S(A_2)\leq b_1+\ldots+ b_{n-2}<b_{n-1}$. This implies that $n-1 \in A_2$. As a consequence, $n-2$ must be in $A_1$, for otherwise we would have $S(A_1)\leq b_n+b_1+b_2+\ldots +  b_{n-3}<2(b_1+b_2+\ldots + b_{n-3})<b_{n-1}\leq S(A_2)$. So $A_1$ contains both $n$ and $n-2$, while $A_2$ contains $n-1$, and we have
\begin{equation*} \sum_{i=0}^{ i<(n-1)/2-1} 2^{2i}+2^{n-3}+\sum_{i\in A_1\setminus \{n,n-2\}} 2^{i-1}=2^{n-2}+\sum_{j\in A_2\setminus\{n-1\}}2^{j-1}\,.
\end{equation*}
Defining now $A'_1=A_1\setminus \{n,n-2\}\cup \{n-1\}$ and  $A'_2=A_2\setminus\{n-1\}\cup\{n-2\}$, again these two sets give a valid counterexample in $\widetilde{\Sigma}_{n-1}$, contradicting the minimality of $n$.

Consider now the case of odd $n$. In this case we can rewrite \eqref{eq:general_n} as
\begin{equation*}2^{n-3} + \sum_{i=0}^{ i<(n-1)/2-1} 2^{2i}+\sum_{i\in A_1\setminus \{n\}} 2^{i-1}=\sum_{j\in A_2}2^{j-1}\,.
\end{equation*}

We notice that $A_2$ must contain either $n-2$ or $n-1$, but not both, because $b_1+b_2+\ldots +b_{n-3}<b_{n}$ but at the same time $b_1+b_2+\ldots + b_{n-3}+b_n<b_{n-2}+b_{n-1}$. Also note that $n-1$ cannot be in $A_1$, for the same reason mentioned in the case of even $n$. So, we are left with the following cases to consider:

a) $n-2\in A_2$, $n-1\notin A_1\cup A_2$ and
\begin{equation*} 2^{n-3} + \sum_{i=0}^{ i<(n-1)/2-1} 2^{2i}+\sum_{i\in A_1\setminus \{n\}} 2^{i-1}=\sum_{j\in A_2\setminus \{n-2\}}2^{j-1} + 2^{n-3}\,,
\end{equation*}
 In this case, by defining $A'_2=A_2\setminus\{n-2\}$ $A'_1=A_1\setminus\{n\}\cup\{n-1\}$ we see that these two sets of indices satisfy $|A'_1|+|A'_2|<(n-1)/2$ and give a counterexample in $\widetilde{\Sigma}_{n-1}$, which contradicts the minimality of $n$.\\

b) $n-2 \in A_1$, $n-1\in A_2$ and
\begin{equation*}2\cdot 2^{n-3} + \sum_{i=0}^{ i<(n-1)/2-1} 2^{2i} +\sum_{i\in A_1\setminus \{n, n-2\}} 2^{i-1}=\sum_{j\in A_2\setminus \{n-1\}}2^{j-1} + 2^{n-2} \,,
\end{equation*}
Here we obtain a counterexample valid for $\widetilde{\Sigma}_{n-1}$ by setting $A'_1=A_1\setminus \{n, n-2\}\cup\{n-1\}$ and $A'_2=A_2\setminus \{n-1\}$.

c) $n-2 \notin A_1\cup A_2$, $n-1\in A_2$ and
\begin{equation*} 2^{n-3} + \sum_{i=0}^{ i< (n-1)/2-1} 2^{2i} +\sum_{i\in A_1\setminus \{n\}} 2^{i-1}=\sum_{j\in A_2\setminus \{n-1\}}2^{j-1} + 2^{n-2} \,.
\end{equation*}
In this case we note that $n-3$ must be in $A_1$, for otherwise $S(A_1)\leq b_n+b_1+\ldots + b_{n-4}<b_n+b_{n-3}<b_{n-1}\leq S(A_2)$ (see Figure \ref{fig:table_b_n}). We can then define $A'_1=A'_1\setminus\{n,n-3\}\cup\{n-1\}$ and $A'_2=A_2\setminus\{n-1\}\cup\{n-3\}$ and again obtain a valid counterexample in $\widetilde{\Sigma}_{n-1}$ which contradicts the minimality of $n$.
\endproof

\begin{Remark}
We note that the condition $|A_1| + |A_2| < n/2$ in the statement of Lemma \ref{Insertion1} is tight, when $n$ is even and greater than or equal to $6$,  because if we take $A_1 = \{b_n\}$ and $A_2 = \{b_{2i+1}:  i=0, \ldots, n/2 - 2\}$  then,  clearly,  $|A_1|+|A_2| = n/2$ and ${S(A_1) = S(A_2)}$.
\end{Remark}

The following corollary follows.
\begin{Corollary}
If $\lambda< 1/4$, $\widetilde{\Sigma}_n$ is $\mathcal{F}_{\lambda,n}$-sum distinct.
\end{Corollary}
The meaning of this Corollary is that it is possible to add one more element to the sequence of powers of two in such a way that it remains $\mathcal{F}_{\lambda,n}$-sum distinct. With the same procedure we can also prove the following statement:
\begin{Lemma}\label{Insertion2}
Let us consider the sequence $\widetilde{\Sigma}_{n}=(b_1,\dots,b_n)$ where, as in Lemma \ref{Insertion1}, we have that
$$b_i:=\begin{cases}
2^{i-1} & i=1,2\dots n-1\\
\sum_{j=0}^{j<n/2-1} 2^{2j} & i=n;
\end{cases}$$
Then, given two subsets $A_1,A_2$ of $[1,n]$ such that $|A_1|+|A_2|< (n-1)/2$,
$$S(A_1)=\sum_{i\in A_1} b_i\not=\sum_{j\in A_2} b_j+2^{n-1}=S(A_2)+2^{n-1}.$$
\end{Lemma}
\proof
We note that the set $\widetilde{\Sigma}_{n}\cup \{2^{n-1}\}$ is $\widetilde{\Sigma}_{n+1}$ whenever $n$ is odd. Therefore, in this case, a contradiction to the statements leads to a contradiction to Lemma \ref{Insertion1} and we can assume $n$ to be even. 

Hence, we suppose now we have a counterexample with $n$ even. We would have that $b_n=\sum_{i=0,\ i\equiv 0\pmod{2}}^{n-4} 2^{i} $ and $n$ must belong to $A_1$. It follows that
$$\sum_{i=0,\ i\equiv 0\pmod{2}}^{n-4} 2^{i} +\sum_{i\in A_1\setminus\{n\}} b_i=\sum_{j\in A_2} b_j+2^{n-1}.$$
Here we note that, since 
$$2^{n-1}=2^{n-2}+2^{n-2}>b_n+\sum_{i=0}^{n-3} 2^{i}=b_n+\sum_{i=1}^{n-2} b_{i},$$
$n-1$ must also belong to $A_1$. 
In this case we would have that:
$$b_n +2^{n-2}+\sum_{i\in A_1\setminus\{n,n-1\}} b_i=\sum_{j\in A_2} b_j+2^{n-1}.$$
We remark that the $(n+1)$-th element of the sequence $\widetilde{\Sigma}_{n+1}$ is $\sum_{i=0,\ i\equiv 0\pmod{2}}^{n-2} 2^{i}$ that is $b_n+2^{n-2}$.
It follows that the set $(\widetilde{\Sigma}_{n}\setminus \{b_n\})\cup \{2^{n-1}, b_n+2^{n-2}\}$ is $\widetilde{\Sigma}_{n+1}$.
Therefore, also for $n$ even, we would obtain a contradiction to Lemma \ref{Insertion1} and thus the statement is verified. 
\endproof

The ideas of Lemmas \ref{Insertion1} and \ref{Insertion2} can be adapted to the following sum-distinct sequence.

\begin{Theorem}[\cite{Lunnon}]\label{Lunnon}
Given $n\geq 67$, there exists a sum-distinct sequence $\overline{\Sigma_{n}}$ of integers such that $c_{1,n}<c_{2,n}<\dots<c_{n,n}$, $0,22\cdot 2^n<c_{n,n}<0,22096\cdot 2^n$ and such that, denoted by $(\overline{c_1},\overline{c_2},\dots,\overline{c_{67}})$ the sequence for $n=67$, we have
$$ c_{i,n}=\begin{cases} 2^{i-1} \mbox{ if } i\leq n-67;\\
2^{n-67}\cdot \overline{c_{i-(n-67)}}\mbox{ otherwise.}
\end{cases}$$
\end{Theorem}
Now we show that it is possible to add one more element also to this sequence in such a way that it remains $\mathcal{F}_{\lambda,n}$-sum distinct. We observe that the proof of our result does not depend on the specific $\bar{c}_i$ values which appear in Theorem \ref{Lunnon}. Indeed, it suffices to analyze only less significant bits in the binary representation of the $c_{i,n}$'s, which are all zeros for $i\geq n-66$.

\begin{figure}
\centering
\includegraphics[scale=1]{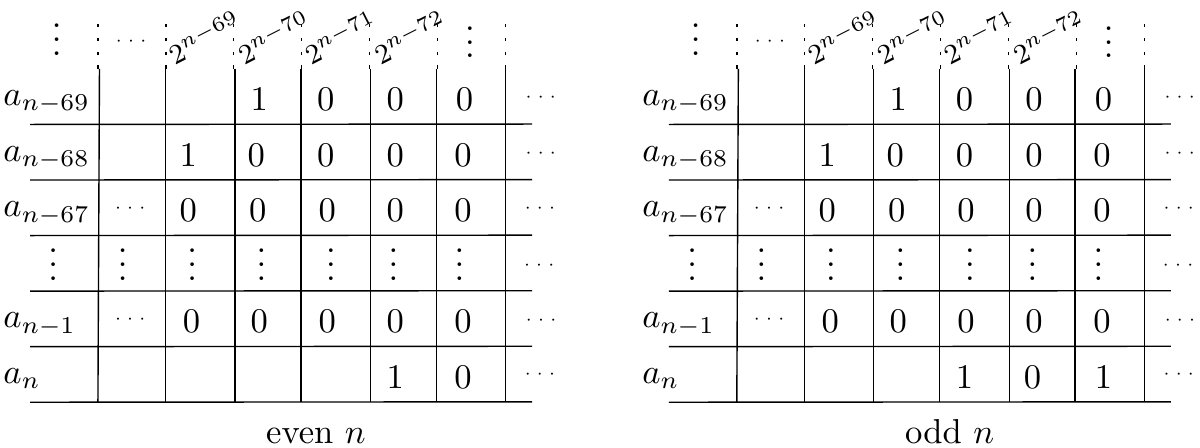}
\caption{Binary representation of the $a_n$ integers used in Proposition \ref{Direct1}.}
\label{fig:table_a_n}
\end{figure}

\begin{Proposition}\label{Direct1}
Let $\Sigma=(a_1,\dots,a_n)$ be the sequence of integers defined by
$$a_i:=\begin{cases}
c_{i,n-1} & i=1\ldots, n-1\\
\sum_{j=0}^{j<(n-68)/2-1} 2^{2j} & i=n
\end{cases}$$
Then, if $\lambda<1/4$ and $n$ is large enough, $\Sigma$ is $\mathcal{F}_{\lambda,n}$-sum distinct.
\end{Proposition}

The structure of the set $\Sigma$ in Proposition \ref{Direct1} is better understood by writing a table of the binary representations of the integers $a_n$, as shown in Figure \ref{fig:table_a_n}.

\proof
Let us suppose, by contradiction, that there exists two disjoint sets $A_1$ and $A_2$ in $\mathcal{F}_{\lambda,n}$ such that $S(A_1)=S(A_2)$.
We note that, if $n\not \in A_1\cup A_2$, we would have two distinct sets of elements of $\overline{\Sigma_{n-1}}$ with the same sums which is in contradiction with the fact that, due to Theorem \ref{Lunnon}, $\overline{\Sigma_{n-1}}$ is sum-distinct.
Therefore, we may assume, without loss of generality, that $n\in A_1$. 
It follows that
$$a_n+\sum_{i\in A_1\setminus \{n\}} a_i=\sum_{j\in A_2}a_j.$$
Set $n':=n-68$.  As a generalization of the method used in Lemma \ref{Insertion1}, we first look at the equation modulo some appropriate power of $2$, namely $2^{n'-1}$ in this case, and then consider possible reminders in the binary expressions for the sums. We set then $A_1':= (A_1\setminus [n',n]) \cup \{n'\}$, $A_2':=A_2\setminus [n',n]$,  and we redefine $a_{n'}$ as $a_{n'}:=\sum_{i=0}^{i<n'/2-1} 2^{2i}$.  Clearly,  both $A_1'$ and $A_2'$ are not empty because $n'\in A_1'$.  Moreover,  since $S(A_1)=S(A_2)$,  $S(A_2') \leq a_1 + a_2 + \ldots + a_{n'-1} < 2^{n'-1}$ and $S(A_1') \leq a_1 + a_2 + \ldots + a_{n'} < 2^{n'}$  we have that either $S(A_1')=S(A_2')$ or $S(A_1')=S(A_2')+2^{n'-1}$. 

In the first case, if $n$ is large enough, we would have that
$$|A_1'|+|A_2'|\leq |A_1| + |A_2| \leq  2\lambda n<n'/2.$$
This would imply that $\widetilde{\Sigma}_{n'}$ is a contradiction to the statement of Lemma \ref{Insertion1} considered for sets $A_1'$, $A_2'$ and for $n'$.

Similarly, in the second case, if $n$ is large enough, we would have that
$$|A_1'|+|A_2'| \leq |A_1| + |A_2| \leq 2\lambda n<(n'-1)/2.$$
Here we would have that $\widetilde{\Sigma}_{n'}$ is a contradiction to the statement of Lemma \ref{Insertion2} considered for sets $A_1'$, $A_2'$ and for $n'$.

Since we obtain a contradiction in all cases, $\Sigma$ is $\mathcal{F}_{\lambda,n}$-sum distinct.
\endproof
In case $\lambda<1/8$ we can even add two elements to the sequence $\overline{\Sigma}$ dividing again the coefficient by $2$. 
At this purpose we need another technical lemma.
\begin{Lemma}\label{sommabinaria}
Let us consider the sequence $(d_1,\dots,d_n)$ where $d_i:=2^{i-1}.$

Then, given three subsets $A_1,A_2,A_3$ of $[1,n]$ such that
$S(A_1)+S(A_2)=S(A_3)$ we have that
$$|A_1|+|A_2|\geq|A_3|.$$
\end{Lemma}
\proof
Assume $A_1$, $A_2$ and $A_3$ form a counterexample with minimum possible value of $|A_1|+|A_2|$. By the uniqueness of the binary representation, it is clear that $Y:= A_1 \cap A_2 \neq \emptyset$.  Also, we note that $n \notin Y$. Then we have
\begin{align*}
\sum_{i\in A_1} 2^{i-1}+\sum_{i\in A_2} 2^{i-1} & = \sum_{i\in A_1\setminus Y } 2^{i-1}+\sum_{i\in A_2\setminus Y} 2^{i-1}+2\sum_{i\in Y} 2^{i-1}\\
& =\sum_{i\in A'_1 } 2^{i-1}+\sum_{i\in A_2'} 2^{i-1}
\end{align*}
where $A'_1=A_1\cup A_2 \setminus Y$ and $A_2'=Y+1$ is obtained by adding $1$ to each element of $Y$. But here $|A_1'|+|A_2'|=|A_1\cup A_2|<|A_1|+|A_2|$, contradicting the assumption that the chosen counterexample minimizes $|A_1|+|A_2|$.
\endproof

\begin{Proposition}\label{prop:Direct2}
Let $n$ be a positive integer, let us set $n'=\lfloor(n-69)/2\rfloor$ and let $\Sigma=(a_1,\dots,a_n)$ be the sequence of integers defined by
$$a_i:=\begin{cases} \sum_{j=\lceil n'/2\rceil-1}^{j<(n-69)/2-1} 2^{2j} & \mbox{ if } i=n;\\
\sum_{j=0}^{j<n'/2-1} 2^{2j} & \mbox{ if } i=n-1;\\
c_{i,n-2} & \mbox{ otherwise.}
\end{cases}$$
Then, if $\lambda<1/8$ and $n$ is large enough, $\Sigma$ is $\mathcal{F}_{\lambda,n}$-sum distinct.
\end{Proposition}

\begin{figure}
\centering
\includegraphics[scale=0.9]{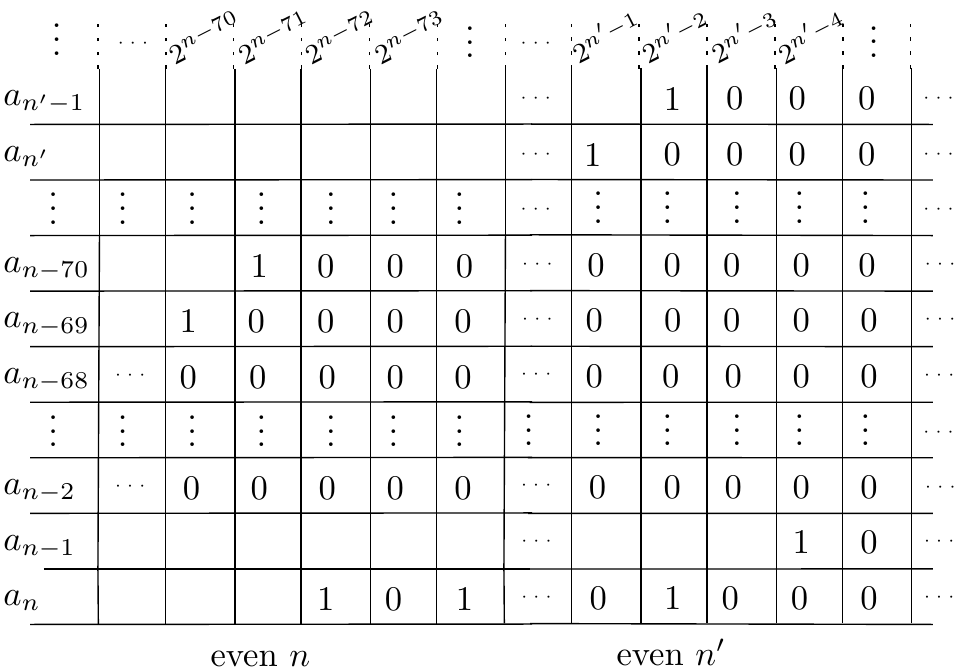}
\caption{Binary representation of the $a_n$ integers used in Proposition \ref{prop:Direct2} when $n$ and $n'$ are even.}
\label{fig:table_a_n_3_11}
\end{figure}

The structure of the set $\Sigma$ in Proposition \ref{prop:Direct2} is better understood by writing a table of the binary representations of the integers $a_n$. In Figure \ref{fig:table_a_n_3_11} we show the table only for even $n$ and $n'$ (the other configurations of $n$ and $n'$ can be easily derived).

\proof
Let us suppose, by contradiction, that there exists $A_1$ and $A_2$ in $\mathcal{F}_{\lambda,n}$ such that $S(A_1)=\sum_{i\in A_1} a_i=\sum_{j\in A_1} a_j=S(A_2)$.
Since $\overline{\Sigma_{n-2}}$ is sum distinct, we may assume, without loss of generality, that $n-1\in A_1$ or $n\in A_1$ and $n-1\not\in A_1, A_2$. Indeed, if both $n$ and $n-1$ do not belong to $A_1\cup A_2$ we would have two distinct sets of elements of $\overline{\Sigma_{n-2}}$ with the same sums which is in contradiction with Theorem \ref{Lunnon}.

In the first case we may assume due to Proposition \ref{Direct1} that $n \notin A_1$ and hence we have $$a_{n-1}+\sum_{i\in A_1\setminus \{n,n-1\}} a_i=\sum_{j\in A_2}a_j.$$

As done in Proposition \ref{Direct1}, we first look at the equation modulo some appropriate power of $2$, namely $2^{n'-2}$ in this case,  and then consider possible reminders in the binary expressions for the sums. 
We set $A_1':=(A_1\setminus [n'-1,n])\cup \{n'\}$, $A_2':=A_2\setminus [n'-1,n]$ and we rename $a_{n'}$ by setting $a_{n'}:=\sum_{i=0}^{i<n'/2-1} 2^{2i}$ where we recall that $n'=\lfloor(n-69)/2\rfloor$. Since $S(A_1)=S(A_2)$,  $S(A_2') \leq a_1 + a_2 + \ldots + a_{n'-2} < 2^{n'-2}$ and $S(A_1') \leq a_1 + a_2 + \ldots + a_{n'-2} + a_{n'} < 2^{n'-1}$ we have that either $S(A_1')=S(A_2')$ or $S(A_1')=S(A_2')+2^{n'-2}$.  In the first case this leads to contradict the statement of Lemma \ref{Insertion1} considered for the sets $A_1'$,  $A_2'$ and for $n'$.  In the second case,  we get a contradiction to the statement of Lemma \ref{Insertion1} considered for the sets $A_1'$,  $A_2' \cup \{ n'-1 \}$ and for $n'$. 

Let us assume now that $n\in A_1$ and $n-1\not\in A_1, A_2$, that is:
\begin{equation}\label{eq:E1}
a_{n}+\sum_{i\in A_1\setminus \{n,n-1\}} a_i=\sum_{j\in A_2}a_j.
\end{equation}
Here we note that by setting $A_3 := \{ 2i + 1 : 0 \leq i < n'/2 -1 \}$,  $a_{n-1} = \sum_{j \in A_3} a_{j}$ so that by adding $a_{n-1}$ to both sides of equation \eqref{eq:E1} we get
\begin{equation}\label{E2} a_{n}+a_{n-1}+\sum_{i\in A_1\setminus \{n,n-1\}} a_i=\sum_{j\in A_2}a_j+\sum_{j \in A_3} a_{j}.\end{equation}

Since $|A_2|<\frac{1}{8}n$,  there exists $h\in [n'-1,n-69]$ that is not in $A_2$ and for which we have that $a_{h} > a_{n-1}$. This implies that 
\begin{equation}\label{E3}
2^{h} > \sum_{\substack{j\in A_2 \\ j <h}} a_j+a_{h} >  \sum_{\substack{j\in A_2\\ j <h}}a_j+ a_{n-1} = \sum_{\substack{j\in A_2\\ j <h}}a_j + \sum_{j \in A_3} a_{j}.
\end{equation} 
Considering the binary representation of the natural numbers  there exists a set $A_2''$ such that
\begin{equation}\label{E4}
\sum_{\substack{j\in A_2, \\ j <h}}a_j + \sum_{j \in A_3} a_{j} = \sum_{j\in A_2''} 2^{j-1}.
\end{equation}
Set $A_1'':= (A_1\setminus \{n\}) \cup \{n-1\}$ and redefine $a_{n-1}$ by setting $a_{n-1} = \sum_{i=0}^{i<(n-69)/2-1} 2^{2i}$. Then thanks to the upper bound of equation \eqref{E3} we know that $A_2'' \subseteq [1, h]$ and hence $a_j = 2^{j-1}$ for $j \in A_2''$ and equation \eqref{E2} can be rewritten as:
\begin{equation}\label{eq:E2}
	\sum_{i\in A_1''} a_i=\sum_{j\in A_2''} a_j + \sum_{\substack{j\in A_2 \\ j > h}} a_j.
\end{equation}
Set $A_2''' := A_2'' \cup (A_2 \setminus [1,h])$.  Since $A_2''$ and $A_2 \setminus [1,h]$ are disjoint,  equation \eqref{eq:E2} becomes
$$
	\sum_{i\in A_1''} a_i=\sum_{j\in A_2'''} a_j.
$$
Now it follows from Lemma \ref{sommabinaria} applied to \eqref{E4} that $|A_2''| \leq |A_2 \setminus [h, n]| + |A_3|$.  Therefore we have that
$$
	|A_2'''| = |A_2''| + |A_2 \setminus [1,h]| \leq |A_2| + |A_3|. 
$$
Moreover, since,  $|A_2|\leq \lambda n<\frac{1}{8}(n-1)$ for $n$ large enough and $|A_3| < \frac{1}{4}(n-1)$,  we obtain that $|A_2'''| <\frac{1}{8}(n-1)+\frac{1}{4}(n-1)$. We also have that, for $n$ large enough, $|A_1''| = |A_1| \leq \lambda n <\frac{1}{8}(n-1)$. Here we note that $|A_1''|+|A_2'''|<\frac{1}{2}(n-1)$ but this is in contradiction with the statement of Proposition \ref{Direct1} considered for the sets $A_1'', A_2'''$ and for $n-1$.
\endproof

As a consequence, we have the following result.
\begin{Theorem}\label{Direct2}
Let $\lambda<1/4$, (resp. $\lambda<1/8$) then, if $n$ is large enough, there exists a sequence $\Sigma=(a_1,\dots,a_n)$ of $\left(\frac{0,22096}{2} \cdot 2^{n}\right)$-bounded integers (resp. $\left(\frac{0,22096}{4} \cdot 2^{n}\right)$-bounded integers) that is $\mathcal{F}_{\lambda,n}$-sum distinct.
\end{Theorem}
\subsection{Multi-Dimensional upper bounds}
In this section we consider the general case $k\geq 1$. First of all, we note that both Theorem \ref{NUB} and Theorem \ref{Direct2} can be used to obtain an upper bound for the $M$ of Problem \ref{SmallSize} also in $\mathbb{Z}^k$.
\begin{Proposition}\label{Recursive}
Let $\bar{\Sigma}$ be an integer $M$-bounded, $\mathcal{F}_{\lambda',n'}$-sum distinct sequence of length $n'$. Then there exists an $M$-bounded sequence $\Sigma$ in $\mathbb{Z}^k$ of length $n$ that is $\mathcal{F}_{\lambda,n}$-sum distinct where $n=kn'$ and $\lambda=\lambda'/k$.
\end{Proposition}

\proof
We set $\bar{\Sigma}_{j}$ to be the sequence in $\mathbb{Z}^k$ whose $j$-th projection is $\bar{\Sigma}$ and that is zero on the other coordinates.
It suffices to consider the sequence $\Sigma=(\bar{\Sigma}_1,\bar{\Sigma}_2,\dots,\bar{\Sigma}_{k})$. Clearly $\Sigma$ is a sequence in $\mathbb{Z}^k$ of length $n$. It is also easy to see that, the existence of $A_1, A_2$ in $\mathcal{F}_{\lambda,n}$ such that $S(A_1)=S(A_2)$ would imply the existence of $A_1', A_2'$ in $\mathcal{F}_{\lambda',n'}$ such that $S(A_1')=S(A_2')$ for $\bar{\Sigma}$. But, since $\bar{\Sigma}$ is an $\mathcal{F}_{\lambda',n'}$-sum distinct sequence, it follows that $\Sigma$ is $\mathcal{F}_{\lambda,n}$-sum distinct.
\endproof

On the other hand, assuming $k>1$, these results can be improved for several values of $\lambda$ using the probabilistic method (see \cite{Alon2}). If $k=1$, instead, the probabilistic method fails to beat the upper bound of Theorem \ref{NUB} (see Remark \ref{rem:Comp}).

We first need another enumerative lemma.
\begin{Lemma}\label{Stima3}
Let $\mathcal{C}$ be the family of the unordered pairs $\{A_1,A_2\}$ of subsets of $[1,n]$ such that:
\begin{itemize}
\item $A_1\cap A_2=\emptyset$;
\item The cardinalities of $A_1$ and $A_2$ are smaller than or equal to $\lambda n$.
\end{itemize}
Then, for $\lambda<1/3$, we have the following upper bound on the cardinality of $\mathcal{C}$

\begin{equation*}\label{eq:multibound2}
|\mathcal{C}|< \frac{\lambda^2 n^2}{2} \cdot 2^{ f(\lambda) n}\,.
\end{equation*}
\end{Lemma}
\proof
It can be easily derived from the proof of Lemma \ref{Stima2}.
\endproof

\begin{Theorem}\label{PUB2}
Let
$$
C_{\lambda, n} = \sqrt[k]{\frac{\lambda^2 n^2}{2 \tau_{\lambda}} 2^{f(\lambda) \tau_{\lambda}}} \text{ and } \tau_{\lambda} = \left \lceil \frac{1}{2^{f(\lambda)}-1} \right \rceil.
$$
Then there exists a sequence $\Sigma=(a_1,\dots,a_n)$, for $n$ large enough,  of $\left( C_{\lambda, n} \cdot 2^{f(\lambda) n/k}\right)$-bounded elements of $\mathbb{Z}^k$ that is $\mathcal{F}_{\lambda,n}$-sum distinct.
\end{Theorem}
\proof
We recall that, if two sets $A_1$ and $A_2$ have the same sum, then also $A_1\setminus (A_1\cap A_2)$ and $A_2\setminus (A_1\cap A_2)$ have the same sum.
Therefore, a sequence $\Sigma$ is $\mathcal{F}_{\lambda,n}$-sum distinct whenever $S(A_1)\not=S(A_2)$ for any $A_1,A_2\in\mathcal{F}_{\lambda,n}$ such that $A_1\cap A_2=\emptyset$. Moreover, since $A_1\not=A_2$, we can assume without loss of generality that $A_2$ is not the empty set.

Now we choose, uniformly at random, the sequence $\Sigma'$ with elements in $[1, M]^k$ and of length $n'$ (whose value will be specified later).  Let $X$ be a random variable that represents the numbers of pairs of elements of $\mathcal{F}_{\lambda,n'}$ such that $A_1\cap A_2=\emptyset$, $S(A_1)=S(A_2)$ and $A_2$ is not the empty set.

Then we need to estimate the following expected value
\begin{align*}
\mathbb{E}[X] & =\mathbb{E}(|\{\{A_1,A_2\}: S(A_1)=S(A_2), A_1,A_2\in \mathcal{F}_{\lambda,n'}, A_1\cap A_2=\emptyset\not= A_2\}|)\\
& = \sum_{\{A_1,A_2\}:\ A_1,A_2\in \mathcal{F}_{\lambda,n'}, A_1\cap A_2=\emptyset\not= A_2} p[S(A_1)=S(A_2)].
\end{align*}
Since $A_1\cap A_2=\emptyset$, the value of $S(A_1)$ is independent from the value of $S(A_2)$. Then the probability $p[S(A_1)=S(A_2)]$ is the following
$$p[S(A_1)=S(A_2)]=\sum_{s\in \mathbb{Z}^k} p[S(A_1)=s]\cdot p[S(A_2)=s].$$
We recall that $A_1\cap A_2=\emptyset\not= A_2$ and hence there exists $i\in A_2\setminus A_1$. Clearly, $A_2$ can sum to $s$ only if $a_i=s-S(A_2\setminus \{i\})$ that happens with probability at most $1/M^k$.
This means that
\begin{align*}
\mathbb{E}[X] & \leq \sum_{\{A_1,A_2\}:\ A_1,A_2\in \mathcal{F}_{\lambda,n}, A_1\cap A_2=\emptyset\not= A_2}\left (\sum_{s\in \mathbb{Z}^k}p[S(A_1)=s] (1/M^k)\right)	\\
& = \frac{1}{M^k}|\{\{A_1,A_2\}:\ A_1,A_2\in \mathcal{F}_{\lambda,n}, A_1\cap A_2=\emptyset\not= A_2\}|.
\end{align*}
Therefore, according to Lemma \ref{Stima3}, we have that
\begin{equation} \label{eq:probabilisticboundE}\mathbb{E}[X]< \frac{1}{M^k} (\lambda n')^2 \cdot 2^{f(\lambda) n' - 1}.\end{equation}
This means that, in case $(1 / M^k) (\lambda n')^2 \cdot 2^{f(\lambda) n' - 1} \leq t$, there exists a sequence $\Sigma'=(a_1,\dots,a_{n'})$ of elements in $\mathbb{Z}^k$ with at most $t$ pairs $\{A_1, A_2\}$ that have the same sum and satisfy the assumptions. Hence, we can remove $t$ elements from $\Sigma'$ and obtain a new sequence $\Sigma=(a_1,\dots,a_n)$, with $n = n' - t$ elements, that is $\mathcal{F}_{\lambda,n}$-sum distinct. Since $n'=n+t$ and due to inequality \eqref{eq:probabilisticboundE}, $\Sigma$ exists whenever
\begin{equation} \label{eq:probabilisticbound}
M \geq (1+o(1)) \sqrt[k]{\frac{\lambda^2 n^2}{2} \frac{2^{f(\lambda) t}}{t}} \cdot 2^{f(\lambda) n/k}.
\end{equation}
It can be seen that the function $g_{\lambda}(t) := \frac{2^{f(\lambda) t}}{t}$ is strictly convex for $t > 0$ and the minimum integer $m$ for which $g_{\lambda}(m+1) \geq g_{\lambda}(m)$ is equal to $\tau_{\lambda}$.  Therefore $t = \tau_{\lambda} = \left \lceil \frac{1}{2^{f(\lambda)}-1} \right \rceil$ is the best choice in order to optimize the inequality \eqref{eq:probabilisticbound}.
\endproof

\begin{Remark}\label{rem:Comp}
We note that for $k=1$ and $n$ sufficiently large the upper bound given in Theorem \ref{NUB} improves the one given in Theorem \ref{PUB2} since 
$$
	 \lambda < \frac{2^{f(\lambda) \tau_{\lambda}}}{2\tau_{\lambda}}\,,
$$
for every $0 < \lambda \leq 1/3$.
\end{Remark}

\begin{figure}%
\centering
\subfloat[\centering $k=1$]{{\includegraphics[width=5.8cm]{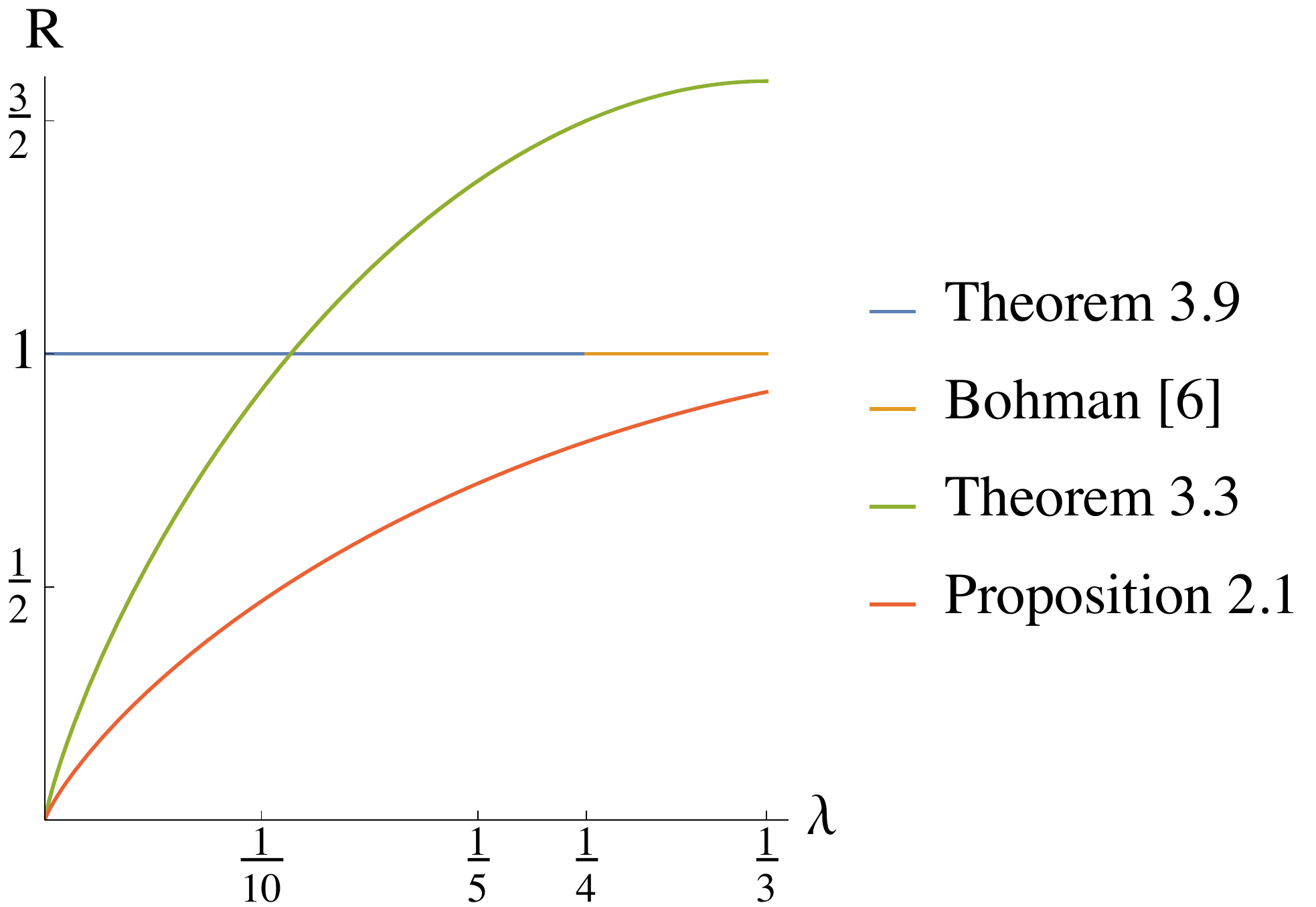} }}%
\qquad
\subfloat[\centering $k>1$]{{\includegraphics[width=5.8cm]{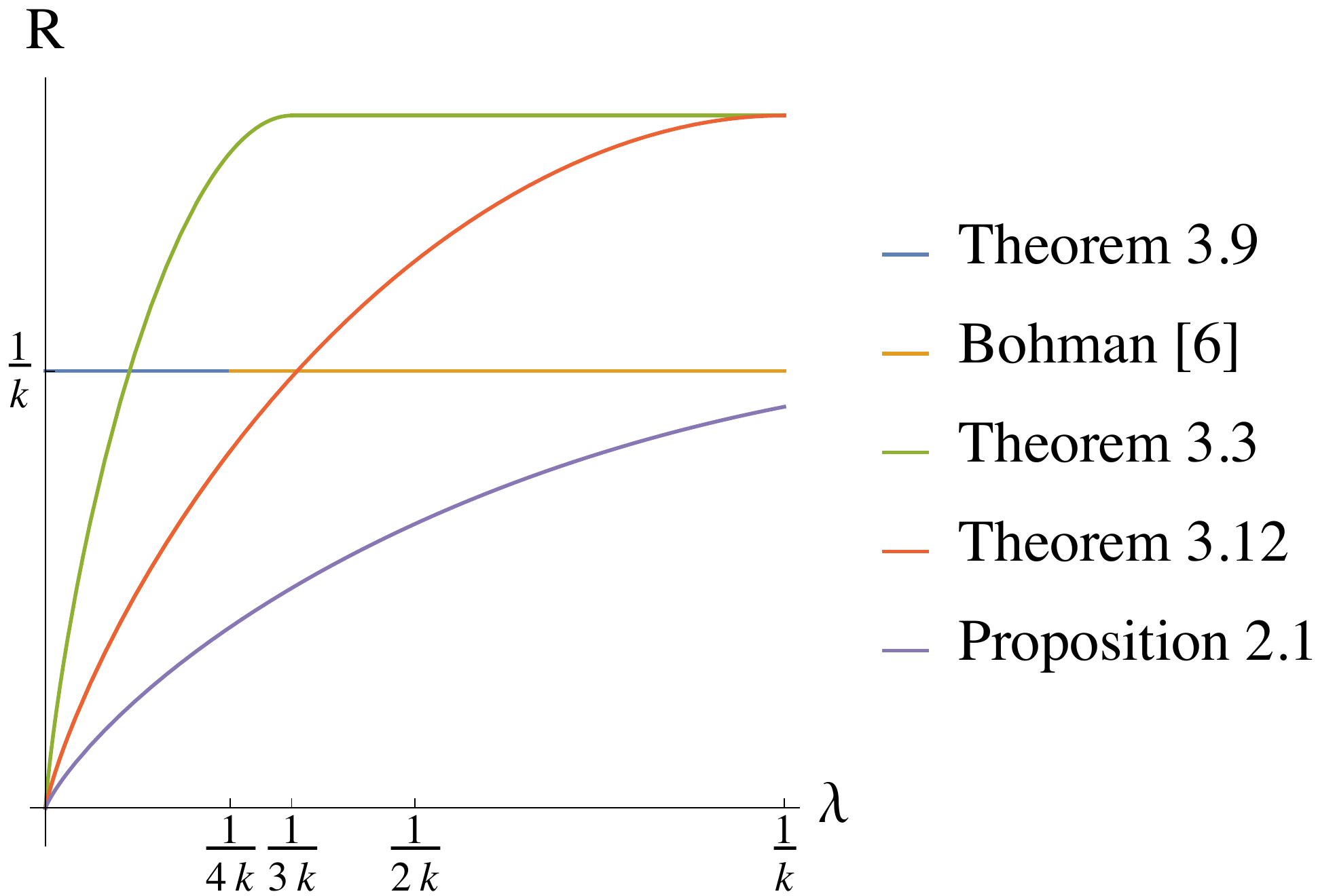} }}%
\caption{Exponent of the upper and lower bounds for $k=1$ and for $k>1$. Here the bounds of Bohman \cite{Bohman}, Theorem \ref{NUB} and Theorem \ref{Direct2} have been extended via Proposition \ref{Recursive}.}%
\label{fig:example2}%
\end{figure}
\section{Conclusions}
In this paper, we investigated a generalization of the Erd\H{o}s sum-distinct problem by weakening the constraint to the family of subsets $\mathcal{F}_{\lambda,n}$ and working in $\mathbb{Z}^k$. We believe that other variations (that are in the same spirit of the problem studied in \cite{MYR}) are also worth considering, such as the following:
\begin{itemize}
\item[1)] $\mathcal{F}$ can be taken as a subfamily of $\mathcal{P}([1,n])$ of a given cardinality,  for example when $\mathcal{F}$ is the family of subsets of $[1,n]$ of  size $n/2$; 
\item[2)] $\mathcal{F}$ can be taken as the family of sets of size at most $m$ (see also \cite{DD} for a similar problem with $m=2$);
\item[3)] each integer is allowed to be covered at most $t$ times by the sums of $\mathcal{F}$.
\end{itemize}
Several of our constructions can be easily adapted to one or more of those situations (for example the probabilistic one works in all those cases). However, no deep idea is needed in this adaptation and, for now, we prefer to keep the treatment more simple and clear. Nevertheless, we plan to investigate those problems more carefully in the future.

\end{document}